\documentclass[12pt,a4paper]{amsart}

\usepackage{hyperref} 
\usepackage{amsfonts}
\usepackage{amsrefs}
\usepackage{amsthm}
\usepackage{amssymb}
\usepackage[all]{xy}

\allowdisplaybreaks

\title[Zero set of semi-invariants for tame canonical algebras]%
  {On the zero set of semi-invariants for regular modules over tame canonical algebras}

\author{Grzegorz Bobi\'nski}

\address{Faculty of Mathematics and Computer Science \\ Nicolaus
Copernicus University \\ ul.~Chopina 12/18 \\ 87-100 Toru\'n \\
Poland}

\email{gregbob@mat.uni.torun.pl}



\newcommand{\bd}{\mathbf{d}}
\newcommand{\be}{\mathbf{e}}
\newcommand{\bh}{\mathbf{h}}
\newcommand{\bm}{\mathbf{m}}
\newcommand{\bx}{\mathbf{x}}
\newcommand{\bP}{\mathbf{P}}
\newcommand{\bQ}{\mathbf{Q}}
\newcommand{\bR}{\mathbf{R}}
\newcommand{\blambda}{\mbox{\boldmath $\lambda$}}

\newcommand{\bbA}{\mathbb{A}}\newcommand{\bbM}{\mathbb{M}}
\newcommand{\bbN}{\mathbb{N}}
\newcommand{\bbP}{\mathbb{P}}
\newcommand{\bbZ}{\mathbb{Z}}

\newcommand{\calA}{\mathcal{A}}
\newcommand{\calC}{\mathcal{C}}
\newcommand{\calO}{\mathcal{O}}
\newcommand{\calP}{\mathcal{P}}
\newcommand{\calQ}{\mathcal{Q}}
\newcommand{\calR}{\mathcal{R}}
\newcommand{\calS}{\mathcal{S}}
\newcommand{\calX}{\mathcal{X}}
\newcommand{\calY}{\mathcal{Y}}
\newcommand{\calZ}{\mathcal{Z}}

\newcommand{\frakC}{\mathfrak{C}}

\let\mod=\undefined
\DeclareMathOperator{\ad}{ad} %
\DeclareMathOperator{\GL}{GL} %
\DeclareMathOperator{\id}{id} %
\DeclareMathOperator{\pd}{pd} %
\DeclareMathOperator{\SI}{SI} %
\DeclareMathOperator{\Ext}{Ext} %
\DeclareMathOperator{\Hom}{Hom} %
\DeclareMathOperator{\mod}{mod} %
\DeclareMathOperator{\gldim}{gl.dim} %
\DeclareMathOperator{\bdim}{\mathbf{dim}} %

\newtheorem*{obs}{Observation}
\newtheorem*{coro}{Corollary}
\newtheorem*{lemm}{Lemma}
\newtheorem*{prop}{Proposition}
\newtheorem*{theo}{Theorem}
\newtheorem*{maintheo}{Main Theorem}

\def\ffrac#1#2{{\textstyle\frac{#1}{#2}}}

\begin{document}

\begin{abstract}
We investigate sets of the common zeros of non-con\-stant
semi-invariants for regular modules over canonical algebras. In
particular, we show that if the considered algebra is tame then
for big enough vectors these sets are complete intersections.
\end{abstract}

\maketitle

Throughout the paper $k$ denotes a fixed algebraically closed
field of characteristic $0$. By $\bbN$ and $\bbZ$ we denote the
sets of non-negative integers and integers, respectively.
Additionally, if $i, j \in \bbZ$, then $[i, j] = \{ l \in \bbZ
\mid i \leq l \leq j \}$.

\makeatletter
\def\@secnumfont{\mdseries} 
\makeatother

\section*{Introduction and the main result} %
\makeatletter
\def\@secnumfont{\mdseries} 
\makeatother

With a finite dimensional algebra $\Lambda$ and a dimension vector
$\bd$ we may associate the variety of $\Lambda$-modules of
dimension vector $\bd$ (see~\ref{var_def}). An interesting problem
investigated in the representation theory of finite dimensional
algebras is the study of geometric properties of these varieties
(see for example~\cites{BobSko1999, BobZwa2002, Bon1991, Bon1998,
DeConStr1981, Gei1996, Kra1982, KraPro1979, Pen1991, PenSko1996,
Zwa2000}). In addition to this topic rings of semi-invariants
(see~\ref{subsect_semiinv}) are also studied (see for
example~\cites{DomLen2000, Hap1984, LeBruPro1990, Kra2001,
Rin1980, SkoWey2000}). Recently, investigations of sets of the
common zeros of non-constant semi-invariants were initiated by
Chang and Weyman (\cite{ChaWey2004}) and then continued by
Riedtmann and Zwara (\cites{Rie2004, RieZwa2003, RieZwa2004,
RieZwa2006}). Their investigations concerned situations of quivers
without relations and were based on known results about
semi-invariants in these cases (among others Sato--Kimura
theorem~\cite{SatKim1977}). An inspiration for their research was
an observation that if, for a given dimension vector, the set of
the common zeros of non-constant semi-invariants has a ``good''
codimension then the coordinate ring of the module variety is free
as a module over the ring of semi-invariants.

An important class of algebras are the canonical algebras
introduced by Ringel~\cite{Rin1984}*{3.7} (see~\ref{can_def}).
These algebras play an important role in representation theory
(see for example~\cites{GeiLen1987, Hap2001, LenPen1999,
Sko1996}). Module varieties over canonical algebras were also
studied (\cites{BarSch2001, Bob2005}). One may distinguish a
special class of modules over canonical algebras, called regular
(see~\ref{reg_def}). The rings of semi-invariants for dimension
vectors of regular modules over canonical algebra were described
by Skowro\'nski and Weyman~\cite{SkoWey1999} (they were also
studied independently by Domokos and Lenzing~\cites{DomLen2000,
DomLen2002}). This description allows to investigate sets $\calZ
(\bd)$ of the common zeros of non-constant semi-invariants for the
dimension vectors $\bd$ of regular modules. The first step in this
direction was made by the author in~\cite{Bob2006}.

If $\bd$ is the dimension vector of a regular module, then we have
a canonical decomposition $\bd = p^{\bd} \bh + \bd'$ of $\bd$
(see~\ref{subsectP}), where $\bh$ is the dimension vector with all
coordinates equal to $1$ and $\bd'$ is the dimension vector of a
regular module such that $\bd' - \bh$ is no longer the dimension
vector of a regular module. Recall that an algebra $\Lambda$ is
called tame if for each dimension $d$ indecomposable modules of
dimension $d$ can be parameterized by a finite number of lines
(see for example~\cite{CB1988}*{Definition~6.5} for a precise
formulation).

The following theorem is the main result of the paper.

\begin{maintheo}
If $\Lambda$ is a tame canonical algebra, then there exists $N$
such that $\calZ (\bd)$ is a complete intersection for all
dimension vectors $\bd$ of regular modules such that $p^{\bd} \geq
N$.
\end{maintheo}

Moreover, we show that also in the case of canonical algebras
there is a connection between the codimension of $\calZ (\bd)$ and
freeness of the coordinate ring over the ring of semi-invariants,
for the dimension vector $\bd$ of a regular module.

The paper is organized as follows. In Section~\ref{sect_can} we
recall necessary facts about quivers, their representations, and
canonical algebras. In Section~\ref{sect_modvar} we present basic
properties of module varieties and rings of semi-invariants. In
particular, we give a description of the sets of the common zeros
of non-constant semi-invariants for the dimension vectors of
regular modules over canonical algebras. In
Section~\ref{sect_proof} we use these results to prove Main
Theorem, and in final Section~\ref{sect_free} we present an
interpretation of the main result in terms of freeness of
coordinate rings over rings of semi-invariants.

For background on the representation theory of algebras we refer
to~\cites{AssSimSko2006,AusReiSma1995}. Basic algebraic geometry
used in the article can be found for example in~\cite{Kun}. Author
gratefully acknowledges the support from the Polish Scientific
Grant KBN No.~1 P03A 018 27. The result presented in this paper
was obtained during the research camp in Szklarska Por\c{e}ba
(June 2006).

\section{Preliminaries on quivers and canonical algebras} %
\label{sect_can}

In this section we present basic facts about quivers and their
representations. We also define canonical algebras and review
their representation theory.

\subsection{}
Recall that by a quiver $\Delta$ we mean a finite set $\Delta_0$
of vertices and a finite set $\Delta_1$ of arrows together with
two maps $s, t : \Delta_1 \to \Delta_0$, which assign to an arrow
$\gamma \in \Delta_1$ its starting and terminating vertex,
respectively. By a path of length $m \geq 1$ in $\Delta$ we mean a
sequence $\sigma = \gamma_1 \cdots \gamma_m$ of arrows such that
$s \gamma_i = t \gamma_{i + 1}$ for $i \in [1, m - 1]$. We write
$s \sigma$ and $t \sigma$ for $s \gamma_m$ and $t \gamma_1$,
respectively. For each vertex $x$ of $\Delta$ we introduce a path
$x$ of length $0$ such that $s x = x = t x$. We only consider
quivers without oriented cycles, i.e., we assume that there is no
path $\sigma$ of positive length such that $t \sigma = s \sigma$.

With a quiver $\Delta$ we associate its path algebra $k \Delta$,
which as a $k$-vector space has a basis formed by all paths in
$\Delta$ and whose multiplication is induced by the composition of
paths. By a relation $\rho$ in $\Delta$ we mean a linear
combination of paths of length at least $2$ with common starting
and terminating vertices. The common starting vertex is denoted by
$s \rho$ and the common terminating vertex by $t \rho$. A set $R$
of relations is called minimal if for every $\rho \in R$, $\rho$
does not belong to the ideal $\langle R \setminus \{ \rho \}
\rangle$ of $k \Delta$ generated by $R \setminus \{ \rho \}$. A
pair $(\Delta, R)$ consisting of a quiver $\Delta$ and a minimal
set of relations $R$ is called a bound quiver. If $(\Delta, R)$ is
a bound quiver, then the algebra $k \Delta / \langle R \rangle$ is
called the path algebra of $(\Delta, R)$.

\subsection{}
By a representation of a bound quiver $(\Delta, R)$ we mean a
collection $M = (M_x, M_\alpha)_{x \in \Delta_0, \, \alpha \in
\Delta_1}$ of finite dimensional vector spaces $M_x$, $x \in
\Delta_0$, and linear maps $M_\alpha : M_{s \alpha} \to M_{t
\alpha}$, $\alpha \in \Delta_1$, such that
\[
\sum_{i \in [1, l]} \lambda_i M_{\alpha_{i, 1}} \cdots
M_{\alpha_{i, m_i}} = 0
\]
for each relation $\sum_{i \in [1, l]} \lambda_i \alpha_{i, 1}
\cdots \alpha_{i, m_i} \in R$. The category of representations of
$(\Delta, R)$ is equivalent to the category of $k \Delta / \langle
R \rangle$-modules (see for
example~\cite{AssSimSko2006}*{Theorem~III.1.6}), and we identify
$k \Delta / \langle R \rangle$-modules and representations of
$(\Delta, R)$. For a representation $M$ its dimension vector
$\bdim M \in \bbZ^{\Delta_0}$ is defined by $(\bdim M)_x = \dim_k
M_x$, $x \in \Delta_0$. For a vertex $x \in \Delta_0$ we denote by
$\be_x$ the corresponding canonical basis vector in
$\bbZ^{\Delta_0}$.

\subsection{}
Let $\Lambda$ be the path algebra of a bound quiver $(\Delta, R)$.
Assume in addition that $\gldim \Lambda \leq 2$. We have the
bilinear form $\langle -, - \rangle = \langle -, - \rangle_\Lambda
: \bbZ^{\Delta_0} \times \bbZ^{\Delta_0} \to \bbZ$ given by
\[
\langle \bd', \bd'' \rangle = \sum_{x \in \Delta_0} d_x' d_x'' -
\sum_{\alpha \in \Delta_1} d_{s \alpha}' d_{t \alpha}'' +
\sum_{\rho \in R} d_{s \rho}' d_{t \rho}''.
\]
It is known (see~\cite{Bon1983}*{2.2}), that if $M$ and $N$ are
$\Lambda$-modules, then
\[
\langle \bdim M, \bdim N \rangle = [M, N] - [M, N]^1 + [M, N]^2,
\]
where following Bongartz~\cite{Bon1994} we write
\begin{align*}
[M, N] & = [M, N]_\Lambda = \dim_k \Hom_\Lambda (M, N),
\\ %
[M, N]^1 & = [M, N]_\Lambda^1 = \dim_k \Ext_\Lambda^1 (M, N),
\\ %
\intertext{and} %
[M, N]^2 & = [M, N]_\Lambda^2 = \dim_k \Ext_\Lambda^2 (M, N).
\end{align*}

\subsection{} \label{can_def}
Let $\bm = (m_1, \ldots, m_n)$, $n \geq 3$, be a sequence of
integers greater than $1$ and let $\blambda = (\lambda_3, \ldots,
\lambda_n)$ be a sequence of pairwise distinct non-zero elements
of $k$ with $\lambda_3 = 1$. By definition $\Lambda (\bm,
\blambda)$ is the path algebra of the bound quiver $(\Delta (\bm),
R (\bm, \blambda))$, where $\Delta (\bm)$ is the quiver
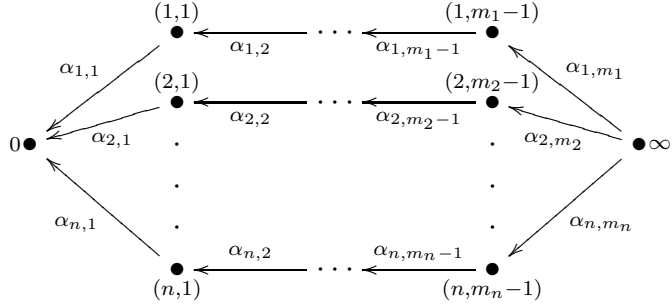
\[
\xymatrix@R=0.25\baselineskip@C=3\baselineskip{%
& \bullet \save*+!D{\scriptstyle (1, 1)} \restore
\ar[lddd]_{\alpha_{1, 1}} & \cdots \ar[l]^-{\alpha_{1, 2}} &
\bullet \save*+!D{\scriptstyle (1, m_1 - 1)} \restore
\ar[l]^-{\alpha_{1, m_1 - 1}}
\\ \\ %
& \bullet \save*+!D{\scriptstyle (2, 1)} \restore
\ar[ld]^{\alpha_{2, 1}} & \cdots \ar[l]^-{\alpha_{2, 2}} & \bullet
\save*+!D{\scriptstyle (2, m_2 - 1)} \restore \ar[l]^-{\alpha_{2,
m_2 - 1}}
\\ %
\bullet \save*+!R{\scriptstyle 0} \restore & \cdot & & \cdot &
\bullet \save*+!L{\scriptstyle \infty} \restore
\ar[luuu]_{\alpha_{1, m_1}} \ar[lu]^{\alpha_{2, m_2}}
\ar[lddd]^{\alpha_{n, m_n}}
\\ %
& \cdot & & \cdot
\\ %
& \cdot & & \cdot
\\ %
& \bullet \save*+!U{\scriptstyle (n, 1)} \restore
\ar[luuu]^{\alpha_{n, 1}} & \cdots \ar[l]_-{\alpha_{n, 2}} &
\bullet \save*+!U{\scriptstyle (n, m_n - 1)} \restore
\ar[l]_-{\alpha_{n, m_n - 1}} }
\]
and $R (\bm, \blambda)$ is the set of the following relations:
\[
\alpha_{1, 1} \cdots \alpha_{1, m_1} + \lambda_i \alpha_{2, 1}
\cdots \alpha_{2, m_2} - \alpha_{i, 1} \cdots \alpha_{i, m_i}, \,
i \in [3, n].
\]
The algebras of the above form are called canonical. In
particular, we say that $\Lambda (\bm, \blambda)$ is a canonical
algebra of type $\bm$. If we fix $\bm$ and $\blambda$, then we
usually write $\Lambda$, $\Delta$, and $R$, instead of $\Lambda
(\bm, \blambda)$, $\Delta (\bm)$, and $R (\bm, \blambda)$,
respectively. From now till the end of the section we assume that
$\Lambda = \Lambda (\bm, \blambda)$ is a fixed canonical algebra.
The following invariant
\[
\delta = \delta_\Lambda = \ffrac{1}{2} \bigl(n - 2 -
\ffrac{1}{m_1} - \cdots - \ffrac{1}{m_n} \bigr)
\]
controls the representation type of $\Lambda$. Namely, $\Lambda$
is tame if and only if $\delta \leq 0$. Moreover, it is known that
$\gldim \Lambda = 2$.

\subsection{}
We abbreviate $\be_{(i, j)}$ by $\be_{i, j}$ for $i \in [1, n]$
and $j \in [1, m_i - 1]$. We put
\[
\bh  = \sum_{x \in \Delta_0} \be_x \qquad \text{and} \qquad
\be_{i, 0} = \bh - (\be_{i, 1} + \cdots + \be_{i, m_i - 1}).
\]
We extend the above definitions by $\be_{i, l m_i + j} = \be_{i,
j}$ for $i \in [1, n]$, $j \in [0, m_i - 1]$, and $l \in \bbZ$.

For $\bd \in \bbZ^{\Delta_0}$ let $\delta_{i, j} (\bd) = d_{i, j -
1} - d_{i, j}$ for $i \in [1, n]$ and $j \in [1, m_i]$. In the
paper we use convention that $d_{i, 0} = d_0$ and $d_{i, m_i} =
d_\infty$ for $\bd \in \bbZ^{\Delta_0}$ and $i \in [1, n]$, and
$d_{i, j} = d_{(i, j)}$ for $i \in [1, n]$ and $j \in [1, m_i -
1]$. Similarly as above we extend this definition by $\delta_{i, l
m_i + j} (\bd) = \delta_{i, j} (\bd)$ for $i \in [1, n]$, $j \in
[1, m_i]$, and $l \in \bbZ$. We also put $\delta_i^{[j_1, j_2]}
(\bd) = \sum_{j \in [j_1, j_2]} \delta_{i, j} (\bd)$ for $i \in
[1, n]$ and $j_1 \leq j_2$. Observe that
\[
\langle \be_{i, j}, \bd \rangle = - \delta_{i, j} (\bd) \qquad
\text{and} \qquad \langle \bd, \be_{i, j} \rangle = \delta_{i, j +
1} (\bd)
\]
for $i \in [1, n]$ and $j \in \bbZ$, and consequently
\[
\langle \be_i^{[j_1, j_2]}, \bd \rangle = - \delta_i^{[j_1, j_2]}
(\bd) \qquad \text{and} \qquad \langle \bd, \be_i^{[j_1, j_2]}
\rangle = \delta_i^{[j_1 + 1, j_2 + 1]} (\bd)
\]
for $i \in [1, n]$ and $j_1 \leq j_2$, where as above
$\be_i^{[j_1, j_2]} = \sum_{j \in [j_1, j_2]} \be_{i, j}$ for $i
\in [1, n]$ and $j_1 \leq j_2$. Finally
\[
\langle \bd, \bh \rangle = d_0 - d_\infty = - \langle \bh, \bd
\rangle.
\]

\subsection{} \label{reg_def}
Let $\calP$ ($\calR$, $\calQ$, respectively) be the subcategory of
all $\Lambda$-modules which are direct sums of indecomposable
$\Lambda$-modules $X$ such that
\[
\langle \bdim X, \bh \rangle > 0 \quad (\langle \bdim X, \bh
\rangle = 0, \, \langle \bdim X, \bh \rangle < 0, \text{
respectively}).
\]
The modules from the category $\calR$ are called regular. We have
the following properties of the above decomposition of the
category of $\Lambda$-modules (see~\cite{Rin1984}*{3.7}).

First, $[N, M] = 0$ and $[M, N]^1 = 0$ if either $N \in \calR \vee
\calQ$ and $M \in \calP$, or $N \in \calQ$ and $M \in \calP \vee
\calR$. Here, for two subcategories $\calX$ and $\calY$ of the
category of $\Lambda$-modules we denote by $\calX \vee \calY$ the
additive closure of their union. Moreover, one knows that
$\pd_\Lambda M \leq 1$ for $M \in \calP \vee \calR$ and
$\id_\Lambda N \leq 1$ for $N \in \calR \vee \calQ$. Secondly,
$\calR$ decomposes into a $\bbP^1 (k)$-family $\coprod_{\lambda
\in \bbP^1 (k)} \calR_\lambda$ of uniserial categories. In
particular, $[M, N] = 0$ and $[M, N]^1 = 0$ if $M \in
\calR_\lambda$ and $N \in \calR_\mu$ for $\lambda \neq \mu$. If
$\lambda \in \bbP^1 (k) \setminus \{ \lambda_1, \ldots, \lambda_n
\}$, where $\lambda_1 = 0$ and $\lambda_2 = \infty$, then there is
a unique (up to isomorphism) simple object $R_\lambda$ in
$\calR_\lambda$ and its dimension vector is $\bh$. On the other
hand, if $\lambda = \lambda_i$ for $i \in [1, n]$, then there are
$m_i$ pairwise non-isomorphic simple objects $R_{i, 0}$, \ldots,
$R_{i, m_i - 1}$ in $\calR_{\lambda_i}$ and their dimension
vectors are $\be_{i, 0}$, \ldots, $\be_{i, m_i - 1}$,
respectively.

For $i \in [1, n]$ and $j_1 \leq j_2$ there is a unique (up to
isomorphism) indecomposable module $R_i^{[j_1, j_2]}$ in
$\calR_{\lambda_i}$ with regular socle (i.e., the socle in the
category $\calR$) $R_{i, j_1}$ and of dimension vector
$\be_i^{[j_1, j_2]}$, where similarly as usual $R_{i, l m_i + j} =
R_{i, j}$ for $i \in [1, n]$, $j \in [0, m_i - 1]$, and $l \in
\bbZ$. Every indecomposable module from $\calR' = \coprod_{i \in
[1, n]} \calR_{\lambda_i}$ is isomorphic to $R_i^{[j_1, j_2]}$ for
some $i \in [1, n]$ and $j_1 \leq j_2$. Moreover, $R_i^{[j_1,
j_2]} \simeq R_i^{[l_1, l_2]}$ if and only if $j_2 - j_1 = l_2 -
l_1$, and $j_1$ and $l_1$ have the same reminder of division by
$m_i$. The regular length (i.e. the length in the category
$\calR$) of $R_i^{[j_1, j_2]}$ is $j_2 - j_1 + 1$ and $\tau
R_i^{[j_1, j_2]} = R_i^{[j_1 - 1, j_2 - 1]}$, where $\tau$ denotes
the Auslander--Reiten translation. We have the following rule of
calculating dimensions of homomorphism spaces between modules in
$\calR'$:
\numberwithin{equation}{subsection} %
\begin{equation} \label{eq_hom}
[R_i^{[j_1, j_2]}, R_i^{[l_1, l_2]}] = \# \{ u \in \bbZ \mid j_1
\leq l_1 + u m_i \leq j_2 \leq l_2 + u m_i \}.
\end{equation}
We also put $\calR'' = \coprod_{\lambda \in \bbP^1 (k) \setminus
\{ \lambda_1, \ldots, \lambda_n \}} \calR_\lambda$.

\subsection{} \label{subsectP}
Let $\bP$, $\bR$ and $\bQ$ denote the sets of the dimension
vectors of modules from $\calP$, $\calR$ and $\calQ$,
respectively. We know from~\cite{Bob2005}*{2.6} that $\bd \in \bP$
($\bd \in \bQ$) if and only if either $\bd = 0$ or $d_0 > d_\infty
\geq 0$ ($0 \leq d_0 < d_\infty$, respectively) and $\delta_{i, j}
(\bd) \geq 0$ ($\delta_{i, j} (\bd) \leq 0$, respectively) for all
$i \in [1, n]$ and $j \in [1, m_i]$.

With a dimension vector $\bd \in \bR$ we may associate its
canonical decomposition (compare~\cite{Rin1980}*{Section~1})
\[
\bd = p^{\bd} \bh + \sum_{i \in [1, n]} \sum_{j \in [0, m_i - 1]}
p_{i, j}^{\bd} \be_{i, j}
\]
in the following way. First, for each $i \in [1, n]$ fix $j_i \in
[0, m_i - 1]$ such that $d_{i, j_i} = \min \{ d_{i, j} \mid j \in
[0, m_i - 1] \}$. Then we put
\begin{align*}
p_{i, j}^{\bd} & =  d_{i, j} - d_{i, j_i}, \; i \in [1, n], \; j
\in [0, m_i - 1],
\\ %
\intertext{and} %
p^{\bd} & = (d_{1, j_1} + \cdots + d_{n, j_n}) - (n - 1) d_0.
\end{align*}
The condition $\bd \in \bR$ implies that $p^{\bd} \geq 0$. We also
put $p_{i, l m_i + j}^{\bd} = p_{i, j}^{\bd}$ for $i \in [1, n]$,
$j \in [0, m_i - 1]$, and $l \in \bbZ$. The canonical
decomposition of $\bd$ is the unique presentation
\[
\bd = p \bh + \sum_{i \in [1, n]} \sum_{j \in [0, m_i - 1]} p_{i,
j} \be_{i, j}
\]
such that $p \geq 0$, $p_{i, j} \geq 0$ for $i \in [1, n]$ and $j
\in [0, m_i - 1]$, and for each $i \in [1, n]$ there exists $j \in
[0, m_i - 1]$ such that $p_{i, j} = 0$.

\makeatletter
\def\@secnumfont{\mdseries} 
\makeatother

\section{Preliminaries on module varieties and semi-invariants}
\label{sect_modvar}

Throughout this section $\Lambda$ is the path algebra of a bound
quiver $(\Delta, R)$.

\subsection{}
\label{var_def} %
For $\bd \in \bbN^{\Delta_0}$ let $\bbA (\bd) = \prod_{\alpha \in
\Delta_1} \bbM (d_{t \alpha}, d_{s \alpha})$. The variety
$\mod_\Lambda (\bd)$ of $\Lambda$-modules of dimension vector
$\bd$ is by definition the subset of $\bbA (\bd)$ formed by all
tuples $(M_\alpha)_{\alpha \in \Delta_1}$ such that
\[
\sum_{i \in [1, l]} \lambda_i M_{\alpha_{i, 1}} \cdots
M_{\alpha_{i, m_i}} = 0
\]
for each relation $\sum_{i \in [1, l]} \lambda_i \alpha_{i, 1}
\cdots \alpha_{i, m_i} \in R$. We identify the points $M$ of
$\mod_\Lambda (\bd)$ with $\Lambda$-modules of dimension vector
$\bd$ by taking $M_x = k^{d_x}$ for $x \in \Delta_0$. The product
$\GL (\bd) = \prod_{x \in \Delta_0} \GL (d_x)$ of general linear
groups acts on $\mod_\Lambda (\bd)$ by conjugations:
\[
(g \cdot M)_\alpha = g_{t \alpha} M_\alpha g_{s \alpha}^{-1}, \;
\alpha \in \Delta_1,
\]
for $g \in \GL (\bd)$ and $M \in \mod_\Lambda (\bd)$. The orbits
with respect to this action correspond bijectively to the
isomorphism classes of $\Lambda$-modules of dimension vector
$\bd$. For $M \in \mod_\Lambda (\bd)$ we denote by $\calO (M)$ the
$\GL (\bd)$-orbit of $M$. It is known (see for
example~\cite{KraRie1986}*{2.2}) that
\[
\dim \calO (M) = \dim \GL (\bd) - [M, M].
\]
We put
\[
a (\bd) = a_\Lambda (\bd) = \dim \bbA (\bd) - \sum_{\rho \in R}
d_{s \rho} d_{t \rho}.
\]
Note that $a (\bd) = \dim \GL (\bd) - \langle \bd, \bd \rangle$
for $\bd \in \bbN^{\Delta_0}$.

\subsection{} \label{subsect_semiinv}
The action of $\GL (\bd)$ on $\mod_\Lambda (\bd)$ induces an
action of $\GL (\bd)$ on the coordinate ring $k [\mod_\Lambda
(\bd)]$ of $\mod_\Lambda (\bd)$ in the usual way:
\[
(g \cdot f) (M) = f (g^{-1} \cdot M)
\]
for $g \in \GL (\bd)$, $f \in k [\mod_\Lambda (\bd)]$, and $M \in
\mod_\Lambda (\bd)$. If $\sigma \in \bbZ^{\Delta_0}$ is a weight,
then we define the weight space
\[
\SI (\Lambda, \bd)_\sigma = \Bigl\{ f \in k [\mod_\Lambda (\bd)]
\mid g \cdot f = \Bigl( \prod_{x \in \Delta_0}
\det\nolimits^{\sigma (x)} (g) \Bigr) f \Bigr\}.
\]
The elements of $\SI (\Lambda, \bd)_\sigma$ are called the
semi-invariants of weight $\sigma$. By the ring of semi-invariants
we mean
\[
\SI (\Lambda, \bd) = \bigoplus_{\substack{\sigma \in
\bbZ^{\Delta_0} \\ \sigma_x = 0 \text{ if } d_x = 0}} \SI
(\Lambda, \bd)_\sigma.
\]
One knows that $\SI (\Lambda, \bd)_0 = k$ (since lack of cycles in
$\Delta$ implies that there is a unique closed orbit in
$\mod_\Lambda (\bd)$). By $\calZ (\bd) = \calZ_\Lambda (\bd)$ we
denote the set of the common zeros of semi-invariants with
non-zero weight for $\bd \in \bR$.

\subsection{}
We present now necessary facts about the rings of semi-invariants
for canonical algebras. For the rest of the section we assume that
$\Lambda = \Lambda (\bm, \blambda)$ is a canonical algebra and
$\Delta = \Delta (\bm)$.

Fix $i \in [1, n]$. An interval $[j_1, j_2]$ with $j_1 < j_2$ is
called $i$-admissible for $\bd \in \bR$ if $p_{i, j_1}^{\bd} =
p_{i, j_2}^{\bd}$ and $p_{i, j}^{\bd} > p_{i, j_1}^{\bd}$ for all
$j \in [j_1 + 1, j_2 - 1]$. Note that $j_2$ is uniquely determined
by $j_1$ and $j_2 \leq j_1 + m_i$. We say that two $i$-admissible
intervals $[j_1, j_2]$ and $[l_1, l_2]$ are equivalent if $j_1$
and $l_1$ have the same reminder of the division by $m_i$
(consequently, $j_2$ and $l_2$ have the same reminder of the
division by $m_i$) --- in other words there exists $u \in \bbZ$
such that $l_1 = j_1 + u m_1$ and $l_2 = j_2 + u m_2$. We will
usually identify equivalent intervals. Let $\calA_i (\bd)$ be the
set of equivalence classes of $i$-admissible intervals for $\bd$
and
\[
\ad (\bd) = \# \calA_1 (\bd) + \cdots + \# \calA_n (\bd).
\]
We will use the following consequence
of~\cite{SkoWey1999}*{Theorem~1.1}.

\begin{prop}
If $\bd \in \bR$, $p^{\bd} \geq n - 1$, and $\mod_\Lambda (\bd)$
is irreducible, then $\SI (\Lambda, \bd)$ is a polynomial ring
generated by $p^{\bd} + 1 + \ad (\bd) - n$ elements.
\end{prop}

If $i \in [1, n]$, $[j_1, j_2] \in \calA_i (\bd)$, and $j \in [0,
m_i - 1]$, then we say that $j$ lies inside $[j_1, j_2]$ if $j_1 +
u m_i \leq j < j_2 + u m_i$ for some $u \in \bbZ$. We will need
the following.

\begin{obs} \label{obs_adm}
Let $\bd \in \bR$, $i \in [1, n]$, and $j \in [0, m_i - 1]$. The
number of $[j_1, j_2] \in \calA_i (\bd)$ such that $j$ lies inside
$[j_1, j_2]$ is bounded above by $p_{i, j + 1}^{\bd} + 1$.
\end{obs}

\begin{proof}
Let $[j_{1, 1}, j_{1, 2}]$, \ldots, $[j_{s, 1}, j_{s, 2}]$ be the
$i$-admissible intervals for $\bd$ with the above property.
Without loss of generality we may assume that
\[
j_{1, 1} < \cdots < j_{s, 1} \leq j < j + 1 \leq j_{s, 2} < \cdots
< j_{1, 2}.
\]
Then $p_{i, j_{s, 2}}^{\bd} > \cdots > p_{i, j_{1, 2}}^{\bd}$ is a
decreasing sequence of $s$ non-negative integers, hence $p_{i, j +
1}^{\bd} \geq p_{i, j_{s, 2}}^{\bd} \geq s - 1$.
\end{proof}

\subsection{}
Now we derive consequences of the connection of semi-invariants
with modules given in~\cite{Dom2002} (see also~\cite{DerWey2002}).
Namely, we have the following description of $\calZ (\bd)$ for
$\bd \in \bR$ with $p^{\bd} > 0$.

\begin{prop}
Let $\bd \in \bR$ and $p^{\bd} > 0$. If $M \in \mod_\Lambda
(\bd)$, then $M \in \calZ (\bd)$ if and only if the following
conditions are satisfied:
\begin{enumerate}

\item
$[R_\lambda, M] \neq 0$ for all $\lambda \neq \lambda_1, \ldots,
\lambda_n$.

\item
$[R_i^{[j_1 + 1, j_2]}, M] \neq 0$ for all $i \in [1, n]$ and
$[j_1, j_2] \in \calA_i (\bd)$.

\end{enumerate}
\end{prop}

By an easy application of the Auslander--Reiten formula
(\cite{AssSimSko2006}*{Theorem~IV.2.13}) we get the following dual
version of the above conditions (see
also~\cite{Dom2002}*{Section~4}).

\begin{obs} \label{obsZd}
Let $\bd \in \bR$ and $M \in \mod_\Lambda (\bd)$.
\begin{enumerate}

\item
If $\lambda \neq \lambda_1, \ldots, \lambda_n$, then
\[
[R_\lambda, M] \neq 0 \qquad {\Leftarrow \joinrel \Rightarrow}
\qquad [M, R_\lambda] \neq 0.
\]

\item
If $i \in [1, n]$, $j_1 < j_2$, and $\delta_i^{[j_1 + 1, j_2]}
(\bd) = 0$, then
\[
[R_i^{[j_1 + 1, j_2]}, M] \neq 0 \qquad {\Leftarrow \joinrel
\Rightarrow} \qquad [M, R_i^{[j_1, j_2 - 1]}] \neq 0.
\]

\end{enumerate}
\end{obs}

Note that for $i \in [1, n]$ and $j_1 < j_2$ the condition
$\delta_i^{[j_1 + 1, j_2]} (\bd) = 0$ is equivalent to $p_{i,
j_1}^{\bd} = p_{i, j_2}^{\bd}$.

\subsection{}
For a subcategory $\calX$ of the category of $\Lambda$-modules and
a dimension vector $\bd$ denote by $\calX (\bd)$ the set of $M \in
\mod_\Lambda (\bd)$ such that $M \in \calX$.

For $\bd \in \bR$ let $\frakC = \frakC (\bd)$ be the set of
quadruples $(\bd', \bd'', [X], q)$ such that $\bd' \in \bP$,
$\bd'' \in \bQ$, $X \in \calR'$, $q \in \bbN$, and $\bd' + \bd'' +
\bdim X + q \bd = \bd$. Observe that $\frakC$ is a finite set. For
$(\bd', \bd'', [X], q) \in \frakC$ let $\calC (\bd', \bd'', [X],
q)$ be the set of $M \in \mod_\Lambda (\bd)$ which are isomorphic
to modules of the form $M' \oplus M'' \oplus X \oplus Y$ with $M'
\in \calP$, $\bdim M' = \bd'$, $M'' \in \calQ$, $\bdim M'' =
\bd''$, and $Y \in \calR''$, $\bdim Y = q \bh$. Obviously
$\mod_\Lambda (\bd)$ is a finite disjoint union of the sets $\calC
(\bd', \bd'', [X], q)$, $(\bd', \bd'', [X], q) \in \frakC$. We
will need the following properties of these sets.

\begin{lemm}
If $\bd \in \bR$ and $(\bd', \bd'', [X], q) \in \frakC$, then
$\calC (\bd', \bd'', [X], q)$ is an irreducible constructible set
of dimension
\[
a (\bd) + \langle \bd - \bd', \bd - \bd'' \rangle - [X, X].
\]
\end{lemm}

\begin{proof}
Compare the proof of~\cite{Bob2006}*{Lemma~3.5}.
\end{proof}

Let $\frakC' = \frakC' (\bd)$ be the set of all $(\bd', \bd'',
[X], q) \in \frakC$ such that the following conditions are
satisfied:
\begin{enumerate}

\item
$\bd' \neq 0$ (equivalently, $\bd'' \neq 0$),

\item
for each $i \in [1, n]$ and each $i$-admissible interval $[j_1,
j_2]$ either $\delta_i^{[j_1 + 1, j_2]} (\bd') > 0$ or $[X,
R_i^{[j_1, j_2 - 1]}] \neq 0$ (equivalently, either
$\delta_i^{[j_1 + 1, j_2]} (\bd'') < 0$ or $[R_i^{[j_1 +1, j_2]},
X] \neq 0$).

\end{enumerate}
Observe that for $M' \in \calP$ the condition $\delta_i^{[j_1 + 1,
j_2]} (\bdim M') > 0$ is equivalent to $[M', R_i^{[j_1, j_2 - 1]}]
\neq 0$. Similarly, for $M'' \in \calQ$ the condition
$\delta_i^{[j_1 + 1, j_2]} (\bdim M'') < 0$ is equivalent to
$[R_i^{[j_1 + 1, j_2]}, M''] \neq 0$.

Another important property, which follows easily
from~\eqref{obsZd} (compare~\cite{Bob2006}*{Lemma~3.6}) is the
following.

\begin{obs}
Let $\bd \in \bR$ and $p^{\bd} > 0$. If $(\bd', \bd'', [X], q) \in
\frakC$, then
\begin{multline*}
\calC (\bd', \bd'', [X], q) \cap \calZ (\bd) \neq \varnothing
\qquad {\Leftarrow \joinrel \Rightarrow} \qquad (\bd', \bd'', [X],
q) \in \frakC'
\\ %
{\Leftarrow \joinrel \Rightarrow} \qquad \calC (\bd', \bd'', [X],
q) \subset \calZ (\bd).
\end{multline*}
\end{obs}

Recall that if $\mod_\Lambda (\bd)$ is irreducible, then it is a
complete intersection of dimension $a (\bd)$ (see for
example~\cite{Bob2005}). Hence we get the following corollary,
which determines our strategy of the proof.

\begin{coro} \label{coro_ineq}
Let $\bd \in \bR$, $p^{\bd} \geq n - 1$, and assume that
$\mod_\Lambda (\bd)$ is irreducible. Then $\calZ (\bd)$ is a
complete intersection provided
\[
[X, X] - \langle \bd - \bd', \bd - \bd'' \rangle \geq p^{\bd} + 1
+ \ad (\bd) - n
\]
for all $(\bd', \bd'', [X], q) \in \frakC'$.
\end{coro}

\makeatletter
\def\@secnumfont{\mdseries} 
\makeatother

\section{Proof of the main result} \label{sect_proof}

Throughout this section $\Lambda = \Lambda (\bm, \blambda)$ is a
fixed canonical algebra and $\Delta$ is its quiver. Our aim in
this section is to prove Main Theorem.

\subsection{}
The first step in our proof is the following.

\begin{lemm}
If $\bd \in \bR$, $(\bd', \bd'', [X], q) \in \frakC'$, and $q >
0$, then there exists $(\bx', \bx'', [X'], q') \in \frakC'$ such
that
\[
[X, X] - \langle \bd - \bd', \bd - \bd'' \rangle > [X', X'] -
\langle \bd - \bx', \bd - \bx'' \rangle.
\]
\end{lemm}

\begin{proof}
Take $\bx' = \bd' + q \bh$, $\bx'' = \bd''$, $X' = X$, and $q' =
0$.
\end{proof}

Let $\frakC''$ be the set of triples $(\bd', \bd'', [X])$ such
that $(\bd', \bd'', [X], 0) \in \frakC'$. We have the following
consequence of the above lemma and Corollary~\ref{coro_ineq}.

\begin{coro} \label{coro_ineq_prim}
Let $\bd \in \bR$, $p^{\bd} \geq n - 1$, and assume that
$\mod_\Lambda (\bd)$ is irreducible. Then $\calZ (\bd)$ is a
complete intersection provided
\[
[X, X] - \langle \bd - \bd', \bd - \bd'' \rangle \geq p^{\bd} + 1
+ \ad (\bd) - n
\]
for all $(\bd', \bd'', [X]) \in \frakC''$.
\end{coro}

\subsection{} \label{subsect_Xcons}
The second, and the most difficult step, is to prove that we may
assume that $p^{\bdim X} = 0$.

Fix $\bd \in \bR$ and $(\bd', \bd'', [X]) \in \frakC''$ such that
$p^{\bdim X} > 0$. We associate to $(\bd', \bd'', [X])$ a new
triple $(\bx', \bx'', [X'])$ such that $\bx' \in \bP$, $\bx'' \in
\bQ$, and $\bx' + \bx'' + \bdim X' = \bd$, in the following way.
Write $X = \bigoplus_{i \in [1, n]} X_i$ with $X_i \in
\calR_{\lambda_i}$, $i \in [1, n]$. Since $p^{\bdim X} > 0$, there
exists $i \in [1, n]$ such that $p^{\bdim X_i} > 0$. Without loss
of generality we may assume that $p^{\bdim X_1} > 0$. Let $j_0$ be
the minimal $j \in [1, m_1]$ such that $\delta_{1, j} (\bd') > 0$,
let $l_2$ be the minimal $l \geq j_0$ such that $R_1^{[j, l]}$ is
a direct summand of $X$ for some $j \leq j_0$ (this definition
makes sense since $p^{\bdim X_1} > 0$), and let $l_1$ be the
minimal $l$ such that $R_1^{[l, l_2]}$ is a direct summand of $X$.
Note that $l_2 < j_0 + m_1$. Write $X = Y \oplus R_1^{[l_1, l_2]}$
and put $\bx' = \bd' + \be_1^{[j_0, l_2]}$, $\bx'' = \bd''$, and
$X' = Y \oplus R_1^{[l_1, j_0 - 1]}$ (where $R_1^{[l_1, j_0 - 1]}
= 0$ if $l_1 = j_0$).

In the following lemma and the next subsection we use the above
notation.

\begin{lemm} \label{lemm_X2}
In the above situation
\[
[X, X] - \langle \bd - \bd', \bd - \bd'' \rangle > [X', X'] -
\langle \bd - \bx', \bd - \bx'' \rangle.
\]
\end{lemm}

\begin{proof}
A crucial role in the proof is played by the following exact
sequence
\[
0 \to R_1^{[l_1, j_0 - 1]} \to R_1^{[l_1, l_2]} \to R_1^{[j_0,
l_2]} \to 0.
\]
By applying the functor $\Hom_\Lambda (-, X)$ to this sequence we
obtain
\begin{multline*}
[R_1^{[l_1, j_0 - 1]}, X] \leq [R_1^{[l_1, l_2]}, X] -
([R_1^{[j_0, l_2]}, X] - [R_1^{[j_0, l_2]}, X]^1)
\\ %
= [R_1^{[l_1, l_2]}, X] - \langle \be_1^{[j_0, l_2]}, \bdim X
\rangle = [R_1^{[l_1, l_2]}, X] + \delta_1^{[j_0, l_2]} (\bdim X).
\end{multline*}
Moreover, by application of the functor $\Hom_\Lambda (R_1^{[l_1,
j_0 - 1]}, -)$ to this sequence we know that
\[
[R_1^{[l_1, j_0 - 1]}, R_1^{[l_1, j_0 - 1]}] \leq [R_1^{[l_1, j_0
- 1]}, R_1^{[l_1, l_2]}],
\]
and consequently
\begin{multline*}
[R_1^{[l_1, j_0 - 1]}, X'] = [R_1^{[l_1, j_0 - 1]}, R_1^{[l_1, j_0
- 1]}] + [R_1^{[l_1, j_0 - 1]}, Y]
\\ %
\leq [R_1^{[l_1, j_0 - 1]}, R_1^{[l_1, l_2]}] + [R_1^{[l_1, j_0 -
1]}, Y] = [R_1^{[l_1, j_0 - 1]}, X].
\end{multline*}
Finally by applying the functor $\Hom_\Lambda (Y, -)$ to the above
sequence we get
\[
[Y, R_1^{[l_1, j_0 - 1]}] \leq [Y, R_1^{[l_1, l_2]}],
\]
hence
\begin{align*}
[X', X'] & = [R_1^{[l_1, j_0 - 1]}, X'] + [Y, R_1^{[l_1, j_0 -
1]}] + [Y, Y]
\\ %
& \leq [R_1^{[l_1, j_0 - 1]}, X] + [Y, R_1^{[l_1, l_2]}] + [Y, Y]
\\ %
& \leq [R_1^{[l_1, l_2]}, X] + \delta_1^{[j_0, l_2]} (\bdim X) +
[Y, R_1^{[l_1, l_2]}] + [Y, Y]
\\ %
& = [X, X]  + \delta_1^{[j_0, l_2]} (\bdim X).
\end{align*}
On other hand
\begin{multline*}
\langle \bd - \bx', \bd - \bx'' \rangle =
\langle \bd - \bd', \bd - \bd'' \rangle - \langle \be_1^{[j_0,
l_2]}, \bd' + \bdim X \rangle
\\ %
= \langle \bd - \bd', \bd - \bd'' \rangle + \delta_1^{[j_0, l_2]}
(\bd') + \delta_1^{[j_0, l_2]} (\bdim X),
\end{multline*}
hence consequently
\[
[X', X'] - \langle \bd - \bx', \bd - \bx'' \rangle \leq [X, X] -
\langle \bd - \bd', \bd - \bd'' \rangle - \delta_1^{[j_0, l_2]}
(\bd')
\]
what finishes the proof.
\end{proof}

\subsection{}
Now we check when $(\bx', \bx'', [X']) \in \frakC''$. For $i \in
[1, n]$ and $[j_1, j_2] \in \calA_i (\bd)$ we say that the triple
$(\bx', \bx'', [X'])$ satisfies the $(i, [j_1,
j_2])$-con\-di\-tion if either $\delta_i^{[j_1 + 1, j_2]} (\bx') >
0$ or $[X', R_i^{[j_1, j_2 - 1]}] \neq 0$ (equivalently, either
$\delta_i^{[j_1 + 1, j_2]} (\bx'') < 0$ or $[R_i^{[j_1 + 1, j_2]},
X'] \neq 0$). Obviously $(\bx', \bx'', [X']) \in \frakC''$ if and
only if $(\bx', \bx'', [X'])$ satisfies $(i, [j_1,
j_2])$-condition for all $i \in [1, n]$ and $[j_1, j_2] \in
\calA_i (\bd)$.

We call a pair $(i, [j_1, j_2])$ consisting of $i \in [1, n]$ and
$[j_1, j_2] \in \calA_i (\bd)$ critical for $(\bd', \bd'', [X])$
if $i = 1$, (after appropriate choice of a representative) $j_2 =
j_0$ and $j_1 < l_1$, $\delta_1^{[j_1 + 1, j_2]} (\bd') = 1$, $[X,
R_i^{[j_1, j_2 - 1]}] = 0$, and $\delta_1^{[j_1 + 1, j_2]} (\bd'')
= 0$. Observe that there may be at most one critical pair for
$(\bd', \bd'', [X])$.

\begin{lemm} \label{lemm_X3}
If $i \in [1, n]$, $[j_1, j_2] \in \calA_i (\bd)$, and $(i, [j_1,
j_2])$ is not a critical pair for $(\bd', \bd'', [X])$, then
$(\bx', \bx'', [X'])$ satisfies the $(i, [j_1, j_2])$-condition.
\end{lemm}

\begin{proof}
If $i \neq 1$ or $\delta_i^{[j_1 + 1, j_2]} (\bd'') < 0$, then the
claim is obvious. Similarly, the claim follows easily if
$\delta_i^{[j_1 + 1, j_2]} (\bd') > 1$, since $\delta_i^{[j_1 + 1,
j_2]} (\bx') \geq \delta_i^{[j_1 + 1, j_2]} (\bd') - 1$. Hence we
may assume that $i = 1$, $\delta_1^{[j_1 + 1, j_2]} (\bd') \leq
1$, and $\delta_1^{[j_1 + 1, j_2]} (\bd'') = 0$.

After an appropriate choice of a representative we may assume that
$j_0 \leq j_2 < j_0 + m_1$. Consider first the case $j_1 \geq
j_0$. If $j_1 \leq l_2 < j_2$, then $\delta_1^{[j_1 + 1, j_2]}
(\bx') > 0$, hence $(\bx', \bx'', [X'])$ satisfies $(1, [j_1,
j_2])$-condition in this case. On the other hand, if either $j_1,
j_2 > l_2$ or $j_1, j_2 \leq l_2$, then $\delta_1^{[j_1 + 1, j_2]}
(\bx') = \delta_1^{[j_1 + 1, j_2]} (\bd')$, thus the claim will
follow if we show that $[X', R_1^{[j_1, j_2 - 1]}] = [X,
R_1^{[j_1, j_2 - 1]}]$ in this case. In order to prove this
equality it is enough to show that $[R_1^{[l_1, l_2]}, R_1^{[j_1,
j_2 - 1]}] = [R_1^{[l_1, j_0 - 1]}, R_1^{[j_1, j_2 - 1]}]$. By
applying the functor $\Hom_\Lambda (-, R_1^{[j_1, j_2 - 1]})$ to
the short exact sequence
\[
0 \to R_1^{[l_1, j_0 - 1]} \to R_1^{[l_1, l_2]} \to R_1^{[j_0,
l_2]} \to 0,
\]
we get a sequence
\begin{multline*}
0 \to \Hom_\Lambda (R_1^{[j_0, l_2]}, R_1^{[j_1, j_2 - 1]}) \to
\Hom_\Lambda (R_1^{[l_1, l_2]}, R_1^{[j_1, j_2 - 1]})
\\ %
\to \Hom_\Lambda (R_1^{[l_1, j_0 - 1]}, R_1^{[j_1, j_2 - 1]}) \to
\Ext_\Lambda^1 (R_1^{[j_0, l_2]}, R_1^{[j_1, j_2 - 1]}).
\end{multline*}
By using~\eqref{eq_hom} and the Auslander--Reiten formula we
obtain that $[R_1^{[j_0, l_2]}, R_1^{[j_1, j_2 - 1]}] = 0$ and
\[
\Ext_\Lambda^1 (R_1^{[j_0, l_2]}, R_1^{[j_1, j_2 - 1]}) =
\Hom_\Lambda (R_1^{[j_1, j_2 - 1]}, R_1^{[j_0 - 1, l_2 - 1]}) = 0,
\]
hence we get the required equality and finish the proof in this
case.

In the second case, i.e.\ when $l_1 \leq j_1 < j_0$, $[R_1^{[l_1,
j_0 - 1]}, R_1^{[j_1, j_2 - 1]}] \neq 0$ and the claim follows
again.

Finally, assume $j_1 < l_1$. If $j_2 > l_2$, then $\delta_1^{[j_1
+ 1, j_2]} (\bx') = \delta_1^{[j_1 + 1, j_2]} (\bd') > 0$. On the
other hand, if $j_2 = j_0$ then $[X, R_1^{[j_1, j_2 - 1]}] \neq
0$, since the pair $(1, [j_1, j_2])$ is not critical. If $j_1 +
m_1 > l_2$ then $[R_1^{[l_1, l_2]}, R_1^{[j_1, j_2 - 1]}] = 0$ and
we get $[X', R_1^{[j_1, j_2 - 1]}] = [X, R_1^{[j_1, j_2 - 1]}]
\neq 0$, while if $j_1 + m_1 \leq l_2$ then $\delta_1^{[j_1 + 1,
j_2]} (\bx') = \delta_1^{[j_1 + 1 + m_1, j_2 + m_2]} (\bx') > 0$,
hence the claim follows in both situations. It remains to consider
the case $j_0 < j_2 \leq l_2$. We show that this situation cannot
happen and this will finish the proof. Indeed, the conditions
$\delta_1^{[j_1 + 1, j_2]} (\bd') \leq 1$ and $\delta_1^{[j_1 + 1,
j_2]} (\bd'') = 0$, which imply $\delta_1^{[j_0 + 1, j_2]} (\bd')
= 0$ (since $\delta_{1, j_0} (\bd') = 1$) and $\delta_1^{[j_0 + 1,
j_2]} (\bd'') = 0$, together with the inequality $d_{1, j_0} >
d_{1, j_2}$, would mean that $p_{1, j_0}^{\bdim X} > p_{1,
j_2}^{\bdim X}$ if $j_0 < j_2 \leq l_2$. As a consequence, there
would exist a direct summand of $X$ of the form $R_1^{[j, l]}$ for
$j \leq j_0 \leq l < j_2 \leq l_2$ in this case
--- a contradiction to the definition of $l_2$.
\end{proof}

\subsection{}
We use now the results of the two previous subsections to make the
next step in the proof.

\begin{lemm}
If $\bd \in \bR$, $(\bd', \bd'', [X]) \in \frakC''$, and $p^{\bdim
X} > 0$, then there exists $(\bx', \bx'', [X']) \in \frakC''$ such
that $\dim_k X' < \dim_k X$ and
\[
[X, X] - \langle \bd - \bd', \bd - \bd'' \rangle \geq [X', X'] -
\langle \bd - \bx', \bd - \bx'' \rangle.
\]
Moreover, the inequality is strict if $p^{\bdim X'} = 0$.
\end{lemm}

\begin{proof}
Without loss of generality we may assume that $p^{\bdim X_1} > 0$,
where $X = \bigoplus_{i \in [1, n]} X_i$ for $X_i \in
\calR_{\lambda_i}$, $i \in [1, n]$. Suppose first there are no
critical pairs for $(\bd', \bd'', [X])$. Then it follows from
Lemmas~\ref{lemm_X2} and~\ref{lemm_X3}, that the triple $(\bx',
\bx'', [X'])$ obtained from $(\bd', \bd'', [X])$ by applying the
construction described in~\eqref{subsect_Xcons} belongs to
$\frakC''$, $\dim_k X' < \dim_k X$, and
\[
[X, X] - \langle \bd - \bd', \bd - \bd'' \rangle > [X', X'] -
\langle \bd - \bx', \bd - \bx'' \rangle.
\]

Assume now that there exists a critical pair for $(\bd', \bd'',
[X])$. Without loss of generality we may assume this pair is of
the form $(1, [j_1, j_0])$ for $j_1 < l_1$, where $j_0$ and $l_1$
(and also $l_2$) have the same meaning as
in~\eqref{subsect_Xcons}. If $R$ is an indecomposable direct
summand of $X$ of the form $R_1^{[u_1, u_2]}$ for $u_1 \leq j_1
\leq u_2$, then it follows that $u_2 \geq j_0$ (since $[X,
R_1^{[j_1, j_0 - 1]}] = 0$), and consequently $u_2 > l_2$ by the
definition of $l_1$ and $l_2$. In particular this means that
$p_{i, j_1}^{\bdim R} \leq p_{i, j_0}^{\bdim R}$ for each direct
summand $R$ of $X$. Since $p_{1, j_1}^{\bdim X} = p_{1,
j_0}^{\bdim X} - 1$, this implies that if $R$ is a direct summand
of $X$ of the form $R_1^{[u_1, u_2]}$ for $u_1 \leq j_0 \leq u_2 <
j_0 + m_1$ different from $R_1^{[l_1, l_2]}$, then $u_1 \leq j_1$
and $u_2 > l_2$.

Let $v_2$ be the minimal $u_2$ such that $R_1^{[u_1, u_2]}$ is a
direct summand of $X$ for $u_1 \leq j_1 \leq u_2$ and let $v_1$ be
the maximal $u_1$ such that $R_1^{[u_1, v_2]}$ is a direct summand
of $X$. Recall that $v_1 < l_1 \leq l_2 < v_2$. Moreover, the
minimality of $v_2$ implies that $v_2 < l_2 + m_1$. Write $X = Y
\oplus R_1^{[l_1, l_2]} \oplus R_1^{[v_1, v_2]}$, and let $X' = Y
\oplus R_1^{[v_1, l_2]} \oplus R_1^{[l_1, v_2]}$. Our definitions
imply that $Y$ has no direct summands of the form $R_1^{[u_1,
u_2]}$ with either $v_1 < u_1 \leq l_1$ and $l_2 < u_2 \leq v_2$,
or $v_1 \leq u_1 < l_1$ and $l_2 \leq u_2 < v_2$, hence
\begin{gather*}
[Y, R_1^{[v_1, l_2]} \oplus R_1^{[l_1, v_2]}] = [Y, R_1^{[l_1,
l_2]} \oplus R_1^{[v_1, v_2]}]
\\ %
\intertext{and} %
[R_1^{[v_1, l_2]} \oplus R_1^{[l_1, v_2]}, Y] = [R_1^{[l_1, l_2]}
\oplus R_1^{[v_1, v_2]}, Y]
\end{gather*}
In addition, tedious analysis shows that
\begin{multline*}
[R_1^{[v_1, l_2]} \oplus R_1^{[l_1, v_2]}, R_1^{[v_1, l_2]} \oplus
R_1^{[l_1, v_2]}] =
\\ %
[R_1^{[l_1, l_2]} \oplus R_1^{[v_1, v_2]}, R_1^{[l_1, l_2]} \oplus
R_1^{[v_1, v_2]}] + 1
\end{multline*}
(here it is important that $v_2 - l_2 < m_1$), hence it follows
that $[X', X'] = [X, X] + 1$. Observe that $(\bd', \bd'', [X'])
\in \frakC''$ (since $[X', R_i^{[u_1, u_2]}] \geq [X, R_i^{[u_1,
u_2]}]$ for all $i \in [1, n]$ and $u_1 \leq u_2$). Moreover,
there are no critical pairs for $(\bd', \bd'', [X'])$, thus it
follows from Lemmas~\ref{lemm_X2} and~\ref{lemm_X3} that for the
triple $(\bx', \bx'', [X''])$ obtained from $(\bd', \bd'', [X'])$
by applying the construction of~\eqref{subsect_Xcons} we have:
$(\bd', \bd'', [X']) \in \frakC''$, $\dim_k X'' < \dim_k X' =
\dim_k X$, and
\begin{multline*}
[X'', X''] - \langle \bd - \bx', \bd - \bx'' \rangle \leq [X', X']
- \langle \bd - \bd', \bd - \bd'' \rangle - 1
\\ %
= [X, X] - \langle \bd - \bd', \bd - \bd'' \rangle.
\end{multline*}
Since $p^{\bdim X''} > 0$, this finishes the proof.
\end{proof}

Let $\frakC'''$ be the set of triples $(\bd', \bd'', [X]) \in
\frakC''$ such that $p^{\bdim X} = 0$. We have the following
consequence of the above lemma and Corollary~\ref{coro_ineq_prim}.

\begin{coro} \label{coro_ineq_bis}
Let $\bd \in \bR$, $p^{\bd} \geq n - 1$, and assume that
$\mod_\Lambda (\bd)$ is irreducible. Then $\calZ (\bd)$ is a
complete intersection provided
\[
[X, X] - \langle \bd - \bd', \bd - \bd'' \rangle \geq p^{\bd} + 1
+ \ad (\bd) - n
\]
for all $(\bd', \bd'', [X]
) \in \frakC'''$.
\end{coro}

\subsection{}
For $\bd \in \bR$ and $(\bd', \bd'', [X]) \in \frakC'''$, let
\begin{align*}
\ad^{(1)} & = \# \{ (i, [j_1, j_2]) \in [1, n] \times \calA_i
(\bd) \mid \delta_i^{[j_1 + 1, j_2]} (\bd') > 0 \},
\\ %
\ad^{(2)} & = \# \{ (i, [j_1, j_2]) \in [1, n] \times \calA_i
(\bd) \mid
\\ %
& \qquad \qquad \qquad \delta_i^{[j_1 + 1, j_2]} (\bd') = 0, \;
\delta_i^{[j_1 + 1, j_2]} (\bd'') < 0 \},
\\ %
\intertext{and} %
\ad^{(3)} & = \# \{ (i, [j_1, j_2]) \in [1, n] \times \calA_i
(\bd) \mid \\ %
& \qquad \qquad \qquad \delta_i^{[j_1 + 1, j_2]} (\bd') = 0 =
\delta_i^{[j_1 + 1, j_2]} (\bd'')\}.
\end{align*}
Obviously
\[
\ad (\bd) = \ad^{(1)} + \ad^{(2)} + \ad^{(3)}.
\]
The final auxiliary step in the proof is as follows.

\begin{lemm} \label{lemm_ineq_ad}
Let $\bd \in \bR$ and $(\bd', \bd'', [X]) \in \frakC'''$.
\begin{enumerate}

\item \label{point1} %
$- \langle \bd, \bd' \rangle \geq (p^{\bd} - n) (d_0' - d_\infty')
+ \ad^{(1)}$,

\item \label{point2} %
$- \langle \bd'', \bdim X \rangle \geq \ad^{(2)}$,

\item \label{point3} %
$[X, X] \geq \langle \bdim X, \bdim X \rangle + \ad^{(3)}$.

\end{enumerate}
\end{lemm}

Before we present the proof of the above lemma, we show how it
implies Main Theorem. Note that
\[
\langle \bd - \bd', \bd - \bd'' \rangle = - \langle \bd', \bd'
\rangle + \langle \bd, \bd' \rangle + \langle \bd'', \bdim X
\rangle + \langle \bdim X, \bdim X \rangle.
\]
Consequently, we have the following corollary being a consequence
of the above lemma and Corollary~\ref{coro_ineq_bis}.

\begin{coro}
Let $\bd \in \bR$, $p^{\bd} \geq n - 1$, and assume that
$\mod_\Lambda (\bd)$ is irreducible. Then $\calZ (\bd)$ is a
complete intersection provided
\[
(\langle \bd', \bd' \rangle  - 1) + (p^{\bd} - n) (\langle \bd',
\bh \rangle - 1) \geq 0
\]
for all $(\bd', \bd'', [X]) \in \frakC'''$.
\end{coro}

\begin{proof}[Proof of Main Theorem]
By repeating arguments used in~\cite{Bob2006}*{Proofs of
Propositions~4.1 and~4.2} we get that
\[
(\langle \bd', \bd' \rangle  - 1) + (p^{\bd} - n) (\langle \bd',
\bh \rangle - 1) - 1 \geq 0
\]
for $\bd' \in \bP$, $\bd' \neq 0$, if $\delta \leq 0$ and $p^{\bd}
\geq N$, where $N = n$ if $\delta < 0$, $N = n + 1$ if $\delta =
0$. Recall from~\cite{BobSko2002}*{Theorem~1} (compare
also~\cite{Bob2005}*{Theorems~1.1 and~1.3~(2)}) that $\mod_\Lambda
(\bd)$ is irreducible, if $\Lambda$ is a tame canonical algebra,
hence the claim follows from the previous corollary.
\end{proof}

\subsection{}
We prove now points~\eqref{point1} and~\eqref{point2} of
Lemma~\ref{lemm_ineq_ad}.

\begin{proof}[Proof of Lemma~\ref{lemm_ineq_ad}~\eqref{point1}]
Let $s = d_0' - d_\infty'$ and $t = d_\infty'$. For each $i \in
[1, n]$ there exists a sequence $0 \leq l_{i, 1} \leq \cdots \leq
l_{i, s} < m_i$ such that
\[
\bd' = t \bh + \sum_{j \in [1, s]} \be (l_{1, j}, \ldots, l_{n,
j})
\]
where for a sequence $(l_1, \ldots, l_n)$ such that $l_i \in [0,
m_i - 1]$ for $i \in [1, n]$ we put
\[
\be (l_1, \ldots, l_n) = \be_0 + \sum_{i \in [1, n]} \sum_{j \in
[1, l_i]} \be_{i, j}.
\]
Note that for $(l_1, \ldots, l_n)$ as above
\[
\langle \bd, \be (l_1, \ldots, l_n) \rangle = - p^{\bd} - \sum_{i
\in [1, n]} p_{i, l_i + 1}^{\bd},
\]
and for $i \in [1, n]$ and $j_1 < j_2$
\[
\delta_i^{[j_1 + 1, j_2]} (\be (l_1, \ldots, l_n)) =
\begin{cases}
1 & l_i \text{ lies inside } [j_1, j_2], \\ %
0 & \text{otherwise}.
\end{cases}
\]
Since for $i \in [1, n]$ and $j_1 < j_2$,
\[
\delta_i^{[j_1 + 1, j_2]} (\bd) > 0 \qquad {\Leftarrow \joinrel
\Rightarrow} \qquad \exists_{j \in [1, s]} \; \delta_i^{[j_1 + 1,
j_2]} (\be (l_{1, j}, \ldots, l_{n, j})) > 0,
\]
the claim follows from Observation~\ref{obs_adm}.
\end{proof}

\begin{proof}[Proof of Lemma~\ref{lemm_ineq_ad}~\eqref{point2}]
Similarly as above for each $i \in [1, n]$ there exists $0 < l_{i,
1} \leq \cdots \leq l_{i, s} \leq m_i$ such that
\[
\bd'' = t \bh + \sum_{j \in [1, s]} \be' (l_{1, j}, \ldots, l_{n,
j})
\]
for $s = d_\infty'' - d_0''$ and $t = d_0''$, where for a sequence
$(l_1, \ldots, l_n)$ such that $l_i \in [1, m_i]$ for $i \in [1,
n]$ we put
\[
\be' (l_1, \ldots, l_n) = \be_\infty + \sum_{i \in [1, n]} \sum_{j
\in [l_i, m_i - 1]} \be_{i, j}.
\]
We also have
\[
\langle \be' (l_1, \ldots, l_n), \bdim X \rangle = - \sum_{i \in
[1, n]} p_{i, l_i - 1}^{\bdim X}
\]
for $(l_1, \ldots, l_n)$ as above.

Fix $i \in [1, n]$. Let $\calA_i^{(2)} (\bd)$ be the set of $[j_1,
j_2] \in \calA_i (\bd)$ such that $\delta_i^{[j_1 + 1, j_2]}
(\bd') = 0$ and $\delta_i^{[j_1 + 1, j_2]} (\bd'') < 0$. We define
a function $f : \calA_i^{(2)} (\bd) \to [1, s]$ given by
\[
f ([j_1, j_2]) = \min \{ j \in [1, s] \mid l_{i, j} - 1 \text{
lies inside } [j_1, j_2] \}.
\]
Our claim will follow if we show that the inverse image of $j \in
[1, s]$ has at most $p_{i, l_{i, j} - 1}^{\bdim X}$ elements. Fix
$j \in [1, s]$ and let $[j_{1, 1}, j_{1, 2}]$, \ldots, $[j_{s, 1},
j_{s, 2}]$ be the intervals in $\calA_i^{(2)}$ whose image under
$f$ is $j$. We may assume that
\[
j_{1, 1} < \cdots < j_{s, 1} \leq l_{i, j} - 1 < l_{i, j} \leq
j_{s, 2} < \cdots < j_{1, 2}.
\]
Then $p_{i, j_{1, 1}}^{\bd} < \cdots < p_{i, j_{s, 1}}^{\bd}$,
hence
\[
s \leq p_{i, j_{s, 1}}^{\bd} - p_{i, j_{1, 1}}^{\bd} + 1 \leq
p_{i, l_{i, j} - 1}^{\bd} - p_{i, j_{1, 1}}^{\bd} + 1.
\]
The definitions of $\calA_i^{(2)} (\bd)$ and $f$ imply that $d_{i,
j_{1, 1}}' = d_{i, l_{i, j} - 1}'$ and $d_{i, j_{1, 1}}'' = d_{i,
l_{i, j} - 1}''$, hence
\begin{multline*}
p_{i, l_{i, j} - 1}^{\bd} - p_{i, j_{1, 1}}^{\bd} = d_{i, l_{i, j}
- 1} - d_{i, j_{1, 1}} = x_{i, l_{i, j} - 1} - x_{i, j_{1, 1}} =
p_{i, l_{i, j} - 1}^{\bx} - p_{i, j_{1, 1}}^{\bx},
\end{multline*}
where $\bx = \bdim X$, thus in order to finish the proof it
remains to show that $p_{i, j_{1, 1}}^{\bx} > 0$. This follows
since the conditions $p_{i, j_{1, 1}}^{\bd} = p_{i, j_{1,
2}}^{\bd}$, $\delta_i^{[j_{1, 1} + 1, j_{1, 2}]} (\bd') = 0$, and
$\delta_i^{[j_{1, 1} + 1, j_{1, 2}]} (\bd'') < 0$, imply that
$p_{i, j_{1, 1}}^{\bx} > p_{i, j_{1, 2}}^{\bx} \geq 0$.
\end{proof}

\subsection{}
Before we give the proof of the last point of
Lemma~\ref{lemm_ineq_ad} we present some auxiliary facts. For $m
\geq 1$ let $A_m$ be the path algebra of the quiver
\[
\Sigma_m = \xymatrix@R=0.25\baselineskip@C=3\baselineskip{%
\bullet \save*+!D{\scriptstyle 1} \restore & \bullet
\save*+!D{\scriptstyle 2} \restore \ar[l] & \cdots \ar[l] &
\bullet \save*+!D{\scriptstyle m - 1} \restore \ar[l] & \bullet
\ar[l] \save*+!D{\scriptstyle m} \restore}.
\]
For an interval $[j_1, j_2]$ with $1 \leq j_1 < j_2 \leq m$ let
$X^{[j_1, j_2]}$ be the unique indecomposable $A_m$-module of
dimension vector $\sum_{j \in [j_1, j_2]} \be_j$. An interval
$[j_1, j_2]$ with $1 \leq j_1 < j_2 \leq m$ is called admissible
for $\bd \in \bbN^{(\Sigma_m)_0}$ if $d_{j_1} = d_{j_2} > 0 $ and
$d_j > d_{j_1}$ for all $j \in [j_1 + 1, j_2 - 1]$. Let $\calA
(\bd)$ be the set of admissible intervals for $\bd$. The following
is a consequence of~\cite{ChaWey2004}*{Theorem~5.7} and the
description of semi-invariants for $\Sigma_m$ obtained
in~\cite{Abe1984} (see also~\cite{Abe1982}).

\begin{prop}
Let $\bd \in \bbN^{(\Sigma_m)_0}$, $\calA'$ be a subset of $\calA
(\bd)$, and $M \in \mod_{A_m} (\bd)$. If $[X^{[j_1 + 1, j_2]},
M]_{A_m} \neq 0$ for all $[j_1, j_2] \in \calA'$, then
\[
[M, M]_{A_m} \geq \langle \bd, \bd \rangle_{A_m} + \# \calA'.
\]
\end{prop}

\begin{proof}[Proof of Lemma~\ref{lemm_ineq_ad}~\eqref{point3}]
For each $i \in [1, n]$ fix $l_i \in [0, m_i - 1]$ such that
$p_{i, l_i}^{\bdim X} = 0$. For $i \in [1, n]$ let $\calS_i$ be
the full subcategory of $\calR_{\lambda_i}$ formed by the objects
$R_i^{[j_1, j_2]}$ such that $l_i < j_1 \leq j_2 < l_i + m_i$. It
is known that there exists an equivalence $F_i$ between $\calS_i$
and the category of $A_{m_i - 1}$-modules such that $F_i
(R_i^{[j_1, j_2]}) = X^{[j_1 - l_i, j_2 - l_i]}$ for $l_i < j_1
\leq j_2 < l_i + m_i$ (in particular, $\bdim (F_i R)_j = p_{i, l_i
+ j}^{\bdim R}$ for $j \in [1, m_i - 1]$). Write $X = \bigoplus_{i
\in [1, n]} X_i$ with $X_i \in \calS_i$. Obviously
\begin{multline*}
[X, X]_\Lambda = [X_1, X_1]_\Lambda + \cdots + [X_n, X_n]_\Lambda
=
\\ %
[F_1 X_1, F_1 X_1]_{A_{m_1 - 1}} + \cdots + [F_n X_n, F_n
X_n]_{A_{m_n - 1}}.
\end{multline*}

Fix $i \in [1, n]$. Note that for each $[j_1, j_2] \in \calA_i
(\bd)$ with $\delta_i^{[j_1 + 1, j_2]} (\bd') = 0 = \delta_i^{[j_1
+ 1, j_2]} (\bd'')$ we have $p_{i, j_1}^{\bdim X_i} = p_{i,
j_2}^{\bdim X_i}$ and $p_{i, j}^{\bdim X_i} > p_{i, j_1}^{\bdim
X}$ for $j \in [j_1 + 1, j_2 - 1]$. Moreover, $[R_i^{[j_1 + 1,
j_2]}, X_i]_\Lambda \neq 0$. This implies in particular that
$p_{i, j_2}^{\bdim X_i} > 0$. Consequently, $[j_1 - l_i, j_2 -
l_i] \in \calA (\bdim F_i X_i)$ and $[X^{[j_1 + 1 - l_i, j_2 -
l_i]}, F_i X_i]_{A_{m_i - 1}} \neq 0$ (here we assume that $[j_1,
j_2]$ is chosen is such a way that $l_1 < j_1 < j_2 < l_1 + m_i$).
Thus it follows from the above proposition that
\[
[F_i X_i, F_i X_i]_{A_{m_i - 1}} \geq \langle \bdim F_i X_i, \bdim
F_i X_i \rangle_{A_{m_i - 1}} + \ad_i^{(3)},
\]
where
\[
\ad_i^{(3)} = \# \{ ([j_1, j_2]) \in \calA_i (\bd) \mid
\delta_i^{[j_1 + 1, j_2]} (\bd') = 0 = \delta_i^{[j_1 + 1, j_2]}
(\bd'')\}.
\]
Since $\langle \bdim F_i X_i, \bdim F_i X_i \rangle_{A_{m_i - 1}}
= \langle \bdim X_i, \bdim X_i \rangle_\Lambda$,
\[
\langle \bdim X, \bdim X \rangle_\Lambda = \langle \bdim X_1,
\bdim X_1 \rangle_\Lambda + \cdots + \langle \bdim X_n, \bdim X_n
\rangle_\Lambda,
\]
and $\ad^{(3)} = \ad_1^{(3)} + \cdots + \ad_n^{(3)}$, the claim
follows.
\end{proof}

\makeatletter
\def\@secnumfont{\mdseries} 
\makeatother

\section{Application to modules of covariants} \label{sect_free}

Let $\Lambda = \Lambda (\bm, \blambda)$ be a canonical algebra and
$\bd \in \bR$ with $p^{\bd} > 0$. The aim of this section is to
prove the following.

\begin{theo}
If $\mod_\Lambda (\bd)$ is irreducible, $\SI (\Lambda, \bd)$ is a
polynomial ring in $s$ variables and the codimension of $\calZ
(\bd)$ in $\mod_\Lambda (\bd)$ equals $s$, then $k [\mod_\Lambda
(\bd)]$ is free as a $\SI (\Lambda, \bd)$-module.
\end{theo}

Note that the conditions of the above theorem are satisfied in the
situations covered by Main Theorem and~\cite{Bob2006}*{Theorem~3}.

The proof of the above theorem basically repeats arguments
from~\cite{Sch1980}*{Proof of Proposition~17.29}.

\begin{proof}
We introduce a grading in $k [\bbA (\bd)]$ in such a way that
polynomials defining $\mod_\Lambda (\bd)$ are homogeneous with
respect to this grading, and consequently $k [\mod_\Lambda (\bd)]$
is graded (recall that the corresponding scheme is reduced
--- see~\cite{Bob2005}*{(3.3)}). Namely, the degree of
$X_{\alpha_{i, j}}^{u, v}$ is $m / m_i$ for $i \in [1, n]$, $j \in
[1, m_i]$, $u \in [1, d_{i, j - 1}]$ and $v \in [1, d_{i, j}]$,
where $m = m_1 \cdots m_n$. Obviously, we may choose generators
$f_1$, \ldots, $f_s$ of $\SI (\Lambda, \bd)$ which are homogeneous
(in fact, one may easily calculate the degrees of the generators
from~\cite{SkoWey1999}*{Theorem~1.1}). It follows from the proof
of~\cite{ZarSam1960}*{Theorem~VII.25, p.~200}, that $f_1$, \ldots,
$f_s$ can be extended to a homogeneous system of parameters for $k
[\mod_\Lambda (\bd)]$. Since $k [\mod_\Lambda (\bd)]$ is a
Cohen--Macaulay ring, the claim is a consequence of arguments
given in~\cite{HocEag1971}*{p.~1036}
\end{proof}

\end{document}